\documentclass{article}
\usepackage{epsfig}

\usepackage{tikz}
\usepackage{mathtools}
\usepackage{upgreek}

\usepackage{amsmath,amssymb,amscd,amsthm,bm} 
\usepackage{mathdots} 

\usepackage{caption}
\usepackage{subcaption}
 
\title{Heron triangles and the hunt for unicorns}
\author{Andrew N.W. Hone\footnote{School of Mathematics, 
Statistics \& Actuarial Science, University of Kent, Canterbury CT2 7FS, UK.\\ 
E-mail: A.N.W.Hone@kent.ac.uk}
}
\date{}

\newcommand{\beq}{\begin{equation}}  
\newcommand{\eeq}{\end{equation}}  
\newcommand{\bea}{\begin{eqnarray}} 
\newcommand{\eea}{\end{eqnarray}}   
\newcommand{\bear}{\begin{array}}  
\newcommand{\eear}{\end{array}}

\newtheorem{thm}{Theorem}[section]

\theoremstyle{definition}

\newcommand{\Q}{{\mathbb Q}}

\newcommand{\rd}{\mathrm{d}}

\newcommand\al{{\alpha}}
\newcommand\be{{\beta}}

\newcommand\gam{{\gamma}}

\newcommand\boa{{\bf a}}
\newcommand\bob{{\bf b}}
\newcommand\boc{{\bf c}}
\newcommand\bok{{\bf k}}

\begin{document} 

\maketitle
 
\begin{abstract} 
A Heron triangle is one that has all integer side lengths and integer area, 
which takes its name from Heron of Alexandria's area formula.  
From a more relaxed point of view, if  rescaling is allowed, then   
one can define a Heron triangle to be one whose side lengths and area are all rational numbers. 
A perfect triangle is 
a Heron triangle with all three medians being rational. According to a longstanding conjecture,   no such triangle exists, so  
perfect triangles are as rare as unicorns.  

However, if perfect is the enemy of good, then perhaps it is best to insist on only two of the 
medians being rational. Buchholz and Rathbun found   an 
infinite family of Heron triangles with two rational medians, 
which they were able to associate with the set of rational points on 
an elliptic curve $E(\Q)$. 
Here we describe a recently discovered  explicit formula for  
the sides, area and medians of these (almost perfect) triangles, 
expressed in terms of a pair of integer sequences: these are Somos sequences, which 
first became popular thanks to David Gale's column in 
\textit{Mathematical Intelligencer}.

\small  
\noindent \textbf{Acknowledgments:} This research was supported by Fellowship EP/M004333/1  from the Engineering \& Physical Sciences Research Council, UK, with EP/V520718/1 UKRI COVID-19 Grant Extension Allocation, and grant IEC\textbackslash R3\textbackslash 193024 from the Royal Society.
\normalsize
\end{abstract} 

\section{Pythagorean triples} 

\setcounter{equation}{0}

One of the oldest problems in the theory of Diophantine equations is to find right-angled 
triangles with integer side lengths, or equivalently triples of positive integers $(a,b,c)$ such that 
\beq \label{pythag} 
a^2 + b^2 = c^2, 
\eeq 
which are called 
Pythagorean triples. 
Examples were known to the Babylonians in around 1800 B.C. 
Taking positive integers $m>n$ and $\uptau$, 
all such triples  
can be determined from  the formula 
\beq\label{euclid} 
a=  \uptau\, (m^2-n^2), \qquad 
b = 2\uptau\, mn, \qquad 
c=\uptau \,(m^2+n^2), 
\eeq 
which was presented by  
Euclid, but without the arbitrary scale factor $\uptau$. A primitive Pythagorean triple is 
one for which $\gcd (a,b,c)=1$, and (up to switching $a$ and $b$) all primitive triples are obtained from (\ref{euclid}) by 
taking $\uptau=1$, and $m,n$ coprime with at least one of them being even.  

The  formula (\ref{euclid}) 
can be derived directly by starting  from simple congruences 
mod 2 and mod 4, 
but another way to obtain it 
is to consider rational points on an algebraic curve, 
namely the unit circle 
$$ 
x^2+y^2=1, \qquad \mathrm{where} \qquad x=\frac{a}{c}, \quad y=\frac{b}{c}. 
$$ 
For any rational point $(x,y)\in\Q^2$ on this circle, 
distinct from the point $(-1,0)$, 
we form the chord joining them, 
given by the 
line $y=t(x+1)$ with slope $t$. Hence 
we find the rational parametrization of the circle, 
\beq\label{ratcirc} 
x=\frac{1-t^2}{1+t^2}, \qquad y =\frac{2t}{1+t^2},  
\eeq 
related to the usual trigonometric parametrization 
$x=\cos\theta$, $y=\sin\theta$ by the ``t-substitution'' of integral calculus, that is 
$t=\tan\tfrac{\theta}{2}$, and formula (\ref{euclid}) follows from taking rational 
$t=\tfrac{n}{m}$ with $0<t<1$.

\section{Heron triangles and unicorns} 

\setcounter{equation}{0}

For a triangle with sides $(a,b,c)$ and semiperimeter $s$, the area 
formula 
\beq\label{heron} 
\Delta=\sqrt{s(s-a)(s-b)(s-c)}, \qquad 
s=\frac{a+b+c}{2}
\eeq
is attributed to Heron of Alexandria. If the side lengths are 
integers 
and the area $\Delta$ is also an integer, then this is called a Heron triangle. 
Allowing the  freedom to rescale all the sides by the same factor, 
it is convenient to define 
a triangle to be Heron whenever the side lengths and the area are all rational numbers. 

Trivially, any right-angled triangle given by a Pythagorean triple is Heron. More generally, 
dropping a perpendicular from any vertex of a Heron triangle splits it into a pair of right-angled triangles with the same height, either joined back-to-back, or overlapping one another, and it is not hard to see that both triangles must have rational sides, so that (up to rescaling) the Heron triangle is built from a pair of  Pythagorean triples. This construction can be used to derive 
a parametric  formula for Heron triangles,  that is 
\beq\label{brahma} 
a=\frac{p^2+r^2}{p}, \qquad 
b=\frac{q^2+r^2}{q}, \qquad 
c=\frac{ (p+q)\, \vert r^2-pq\vert}{pq},  \qquad \mathrm{with} \quad \Delta=rc, 
\eeq 
for arbitrary positive rational numbers $p,q,r$ such that $r^2\neq pq$,  
which was known to Brahmagupta in the 7th century A.D.\ \cite{dickson}. 

\begin{figure}
 \centering
\epsfig{file=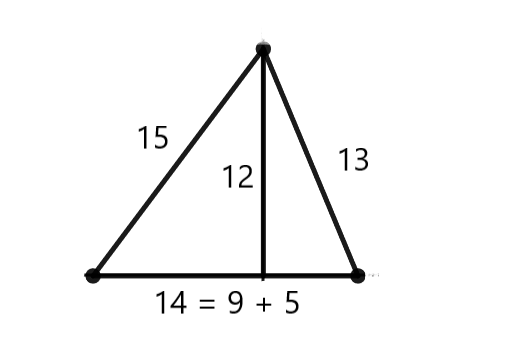, height=1.5in, width=2.5in}
\caption{The (15,13,14) Heron triangle from the two Pythagorean triples 
(9,12,15) \& (5,12,13).} 
\label{151314fig}
\end{figure}

To see a particular example, taking $p=3,q=4,r=6$ in (\ref{brahma}) leads to the Heron triangle with 
$a=15$, $b=13$, $c=14$ and area $\Delta=84$, which can be 
built out of 
the $(5,12,13)$ and $(9,12,15)$ Pythagorean triples, by placing these 
two right-angled triangles back-to-back along the height $2r=12$, as in Fig.\ref{151314fig}. (This choice of parameters is not unique: for instance,  
ordering the sides differently as $(a,b,c)=(15,14,13)$ 
instead of $ (15,13,14)$,
 the same Heron triangle can be obtained from 
$p=\tfrac{147}{13}$, $q=\tfrac{126}{13}$, $r=\tfrac{84}{13}$.)
A systematic method 
for enumerating Heron triangles with integer sides was given by Schubert \cite{schubert}.  
In Schubert's scheme, $(15,13,14)$ is the first example of a Heron triangle with integer sides which is not right-angled or isosceles. However, if we combine the same two Pythagorean triples 
by overlapping the triangles (rather than back-to-back as in the figure), then 
we get the $(15,13,4)$ Heron triangle with a smaller area $\Delta=24$, 
which nevertheless appears further down in Schubert's list.

The unicorns in our story are perfect triangles: triangles which have 
all integer sides, all integer medians, and integer area. Does a perfect triangle exist, or equivalently, is there a Heron triangle 
with three rational medians? It is believed that there is no such thing, 
but despite incorrect ``proofs'' in the literature, the problem remains open \cite{guy}. 
The rest of our discussion is devoted to seeing how close we can get to perfection. 

Henceforth the medians bisecting sides $a,b,c$ are denoted $k,\ell,m$, respectively, 
which leads to the relations 
\beq\label{kmed} 
k^2=\frac{1}{4}(2b^2+2c^2-a^2), \quad  
\ell^2=\frac{1}{4}(2c^2+2a^2-b^2), \quad 
m^2=\frac{1}{4}(2a^2+2b^2-c^2). 
\eeq 
We label the angles adjacent to the median $k$ 
as in Fig.\ref{triangle},  and our first step towards the elusive perfect triangle will be to  
consider the requirement that just this median should be rational. 
 

\section{Heron triangles with one rational median} 

\setcounter{equation}{0}

From a construction of parallelograms with rational sides, area and diagonals, 
Schubert was led to the case of Heron triangles with one  median being rational, 
and went on to present an argument that such triangles could not have a second rational median, which a fortiori would rule out  the existence of perfect triangles. However, 
as pointed out 
by Dickson \cite{dickson}, this argument contained an oversight. This flaw notwithstanding, an identity of Schubert for Heron triangles with one rational median 
is crucial for what follows. 

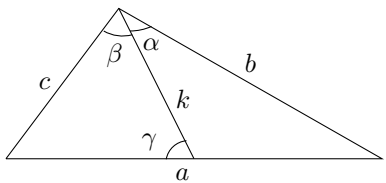
\begin{figure}
\centering
\begin{tikzpicture}
\draw (2,0) 
  -- (7,0) 
-- (5.25,1) node[above=1pt] {$b$}
--(3.5,2) 
--(2.75,1) node[left=1pt] {$c$}
-- (2,0); 
\draw (3.5,2)  
--(4.1,0.8) node[right=1pt] {$k$} 
-- (4.5,0) ;
\draw (3.3,1.7) arc [start angle=240, end angle=285, radius=5mm] node[anchor=north east] {$\be$};
\draw (3.65,1.7) arc [start angle=270, end angle=294, radius=7mm] node[below=1pt] {$\al$};
\draw (4.4,0.24) arc [start angle=95, end angle=170, radius=3mm]  node[anchor=north west] {$a$} node[anchor=south east] {$\gam $};
\end{tikzpicture}
\caption{Triangle with one labelled median}
\label{triangle}
\end{figure}

If we write $\bob,\boc,\bok$ for the vectors corresponding to the lengths $b,c,k$,  
respectively, directed outwards  from the top vertex in Fig.\ref{triangle}, and 
$\boa = \bob-\boc=2(\bok-\boc)=2(\bob-\bok)$, then the dot product 
$(\bob-\boc)\cdot \bok = \boa\cdot\bok$ gives $bk\cos\al-ck\cos\be=ak\cos\gam$, 
while the area of the triangle is $\Delta
=\vert \bob \times \bok\vert 
=\vert \boc \times \bok\vert 
=\tfrac{1}{2}\vert \boa \times \bok\vert$, 
which gives    
$\Delta = bk\sin\al =ck\sin\be =\frac{1}{2}ak\sin\gam$; 
so combining these relations produces the identity 
\beq\label{schubertid} 
2\cot\gam = \cot\al-\cot\be
. 
\eeq 

Given three angles $\al,\be,\gam$ 
in the interval $(0,\pi)$, subject to $\al+\be<\pi$,
it is convenient to take 
\beq\label{mpx} 
M=\mathrm{cot}\tfrac{\al}{2}, \qquad 
P=\mathrm{cot}\tfrac{\be}{2}, 
\qquad 
X=\mathrm{cot}\tfrac{\gam}{2}, 
\eeq
as parameters, 
and then  by standard trigonometric identities 
(equivalent to the ``t-substitution'' in (\ref{ratcirc}) above) 
the identity (\ref{schubertid}) becomes a rational relation between 
these three quantities, namely 
\beq\label{schub}
M-\frac{1}{M}=P-\frac{1}{P}+
2\left( X-\frac{1}{X}\right). 
\eeq
This gives the equation of a surface in 
three-dimensional space with 
coordinates $(M,P,X)$, which can be   
rewritten as the vanishing  of a  polynomial: 
$
2MP(X^2-1)+MX(P^2-1)-PX(M^2-1)=0; 
$
we refer to it 
as the Schubert surface. 

From the  half-angle identity $\mathrm{cot}\tfrac{\al}{2}=\sin\al/(1-\cos\al)$ we have 
$M=\Delta/(bk-\bob\cdot\bok)$. Using  the analogous expressions for $P$ and $X$,  
together with dot product relations, we can express
these Schubert parameters in terms of the area, side lengths and the median 
by 
the formulae 
\beq\label{schubb} 
M=\frac{4\Delta}{4bk+a^2-3b^2-c^2}, \quad 
P=\frac{4\Delta}{4ck+a^2-b^2-3c^2}, \quad 
X=\frac{4\Delta}{2ak-b^2+c^2}. 
\eeq  
The ratios of the side lengths are given in terms of the  Schubert parameters by 
\beq\label{aboverc} 
\frac{a}{c}=\frac{2(X+X^{-1})}{P+P^{-1}}, \qquad
\frac{b}{c}=\frac{M+M^{-1}}{P+P^{-1}}.
\eeq

The formulae (\ref{schubb}) show that each Heron triangle with a rational median $k$ 
produces a rational point  on the Schubert surface (\ref{schub}), 
with positive coordinates $(M,P,X)\in \Q^3$. How about the converse: does every rational  
point on this surface correspond to a Heron triangle with (at least) one rational median? 
In fact, using certain discrete symmetries of the surface (sending $(M,P,X)\to(M^{-1},P^{-1},X^{-1})$, 
or replacing one of the Schubert parameters by minus its reciprocal), we can start with any triple of 
non-zero values
$(M,P,X)\in \Q^3$ satisfying (\ref{schub}), and turn it into a valid positive triple. 
Then 
the side lengths $(a,b,c)$ are determined by $(M,P,X)$ using the rational 
expressions (\ref{aboverc}), up to an 
arbitrary choice of scale; after fixing the scale, any two of the equations (\ref{schubb}) allow 
the rational numbers $k$ and $\Delta$ to be recovered.

\section{Triangles with two rational medians} 

\setcounter{equation}{0}

In striving to get closer to perfection, another possible 
direction for our first step is to drop the initial requirement that the area $\Delta$ should be rational, 
and just consider triangles with rational sides $(a,b,c)$ and two rational medians $k,\ell$.   
In his PhD thesis, 
Buchholz obtained a rational parametrization of all such 
triangles, given by 
the formulae 
\beq\label{abcparam} 
\begin{array}{rcl} 
a & = &\uptau \,(-2\theta^2\phi -\theta\phi^2+2\theta\phi-\phi^2+\theta+1), \\
b & = & \uptau\,(\theta^2\phi+2\theta\phi^2-\theta^2+2\theta\phi-\phi+1), \\
c & = & \uptau\, (\theta^2\phi-\theta\phi^2+\theta^2+2\theta\phi+\phi^2+\theta-\phi),
\end{array} 
\eeq 
where 
$\theta,\phi$ are rational numbers subject to constraints 
ensuring positivity of the side lengths, namely 
\beq\label{thphcon} 
0<\theta<1, \qquad 0<\phi<1, \qquad \phi+2\theta>1, 
\eeq 
and the positive parameter 
$\uptau\in\Q$ allows for the arbitrary choice of scale. 
Conversely, the parameters $(\theta,\phi)\in\Q^2$ 
can be written as rational functions of the  side lengths and two medians, given by 
\beq\label{thetaphi}
\theta = \frac{c-a\pm2\ell}{2s}, \qquad 
\phi   = \frac{b-c\pm2k}{2s},
\eeq 
where $s=\frac{1}{2}(a+b+c)$ is the semiperimeter, as before. 

Note that in (\ref{thetaphi}) there are two independent choices of $\pm$ signs, and hence  four different pairs $(\theta,\phi)$ associated with the same 
rational triangle with two rational medians.

\section{Intermezzo: Somos-5 sequences} 

\setcounter{equation}{0}

Before we continue  on  our quest for the perfect triangle, we must recall some beautiful observations 
made by Michael Somos \cite{somos}. The saga of Somos sequences attracted 
widespread attention due to \textit{Mathematical Intelligencer} articles by David Gale \cite{gale}, and provided inspiration for 
the study of the Laurent phenomenon 
and its development in Fomin and Zelevinksy's theory of cluster algebras   \cite{fz, fordy_marsh},
which has been one of the hottest topics in algebra for almost 25 years.   

A recurrence relation of Somos type is a homogeneous quadratic recurrence relation of a particular form. 
Here we focus on the example of Somos-5, which is the recurrence relation of order 5 given by 
\beq\label{s5recu} 
S_{n+5}S_n = S_{n+4}S_{n+1}+S_{n+3}S_{n+2} 
\eeq 
Somos noticed that if all five initial values are 
1, then the resulting sequence\footnote{
See 
{\tt http://oeis.org/A006721}.}  
begins with  
\beq\label{somos5orig} 
1, 1, 1, 1, 1, 2, 3, 5, 11, 37, 83, 274, 1217, 6161, 22833, 165713, \ldots, 
\eeq 
and consists entirely of integers.  
This seems very surprising, because at each iteration of (\ref{s5recu}) one must divide the 
right-hand side by $S_n$ to obtain the new term $S_{n+5}$. The Laurent property 
provides one explanation for the integrality of the sequence (\ref{somos5orig}): 
if the initial values $S_j$, $1\leq j\leq 5$ for  the recurrence are considered as variables, 
then each iterate turns out to be a  polynomial in these quantities and their reciprocals with 
integer coefficients: 
$S_n ={\mathcal P}_n (S_{1}^{\pm 1},  S_{2}^{\pm 1},S_{3}^{\pm 1},S_{4}^{\pm 1},S_{5}^{\pm 1})$ (that is, a Laurent polynomial). 
Upon substituting 
$S_1=S_2=S_3=S_4=S_5=1$ into each polynomial ${\mathcal P}_n$, the 
integer sequence (\ref{somos5orig}) results. 


\begin{figure}
 \centering
\epsfig{file=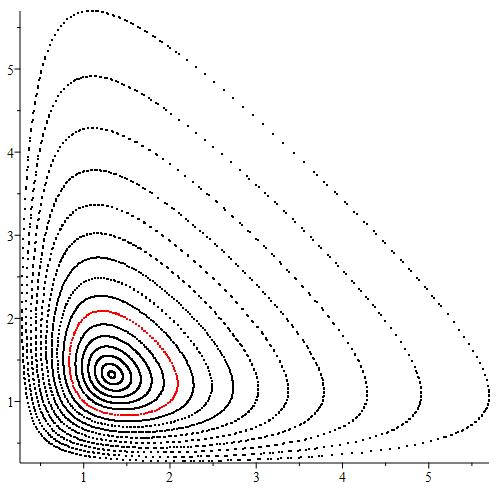, height=2in, width=2in}
\caption{Some orbits of the map (\ref{qrts5}) in the positive quadrant.} 
\label{s5orb}
\end{figure} 

Another completely different way to understand Somos-5 sequences  relies on a connection with integrable maps, 
which are discrete analogues of exactly solvable systems in Hamiltonian mechanics. 
To see this connection, note that the recurrence (\ref{s5recu}) has three independent scaling symmetries: 
rescaling even/odd index terms separately, so $S_{2j}\to A_+\,S_{2j}$,  $S_{2j+1}\to A_-\,S_{2j+1}$, and rescaling 
$S_n\to B^n\,S_n$ for any $n$, where $A_+,A_-,B$ are arbitrary non-zero constants. Moreover, we can form a sequence of ratios 
that is left invariant by these scaling symmetries, 
and find that it satisfies a recurrence of second order: 
\beq\label{ratun}
u_n =  \frac{S_{n-2}S_{n+1}}{S_{n-1}S_{n}} \implies u_{n+1} u_{n-1} = 1+\frac{1}{u_{n}}. 
\eeq  
By considering  $(U,V)=(u_n,u_{n+1})$ as a point in the plane, each shift $n\to n+1$ of the discrete ``time''  in 
(\ref{ratun}) 
is equivalent to an iteration of a birational transformation (a rational map  with a rational inverse):
\beq\label{qrts5} 
\varphi: \qquad \left(\begin{array}{c} U \\ V \end{array} \right) \mapsto 
\left(\begin{array}{c} V \\ U^{-1}(1+ V^{-1}) \end{array} \right)  
.
\eeq 

The transformation (\ref{qrts5}) is an example of a Quispel-Roberts-Thompson (QRT) map: 
such maps have arisen in various physical contexts, including statistical mechanics,  nonlinear waves (solitons) and quantum field theory \cite{qrt}. 
In a suitable regime, the iterates of the map appear like a stroboscopic view of a mechanical system with one degree of freedom. More precisely, 
the transformation $\varphi$ is area-preserving (symplectic): it preserves the logarithmic area element $(UV)^{-1}\,\rd U\,\rd V$ in the plane; 
and it obeys conservation of energy, where ``energy'' in this case is the rational function 
\beq\label{Jtilde} 
\tilde{J}= U+V+\frac{1}{U}+\frac{1}{V}+\frac{1}{UV}. 
\eeq     
The  level sets of this function  are plane curves  $\tilde{J}=\,$constant, and each orbit of $\varphi$ lies on a fixed level set. 
The behaviour is especially  regular in the positive quadrant $U>0$, $V>0$, where each orbit densely fills a compact oval 
(see Fig.\ref{s5orb}, where 300 points are plotted on each orbit). 

\begin{figure}
 \centering
\epsfig{file=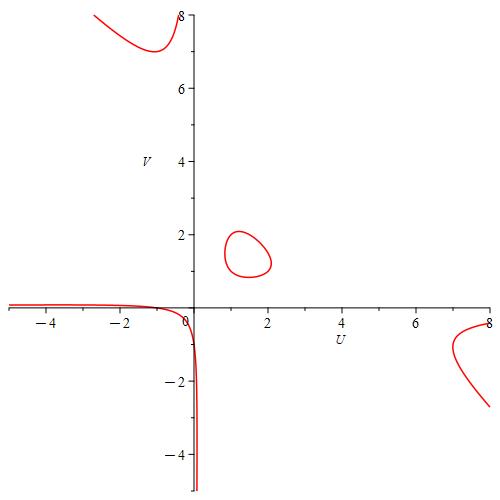, height=2in, width=2in}
\caption{The curve (\ref{uvcurve}) in the $(U,V)$-plane.} 
\label{s5curvep}
\end{figure}

We shall see that in relation to Heron triangles with two rational medians, two different 
integer sequences appear, namely  the pair of    
Somos-5 sequences 
given by 
\beq\label{s5seq} 
(S_n): \quad1,1,1,2,3,5,11,37,83,274,\ldots , \qquad (T_n): \quad 0,1,-1,1,1,-7,8,-1,-57,391,\ldots , 
\eeq 
where the terms above are listed starting from the index $n=0$. 
The first one,  $(S_n)$, is just the original Somos-5 sequence (\ref{somos5orig}), but indexed differently: 
it corresponds to the orbit of the 
map $\varphi$ through the point $(1,1)$, while the second sequence,  $(T_n)$, corresponds to the orbit through the point $(-1,7)$.  It 
is easily verified that 
both of these orbits lie on the same level curve $\tilde{J}=5$ of the function (\ref{Jtilde}),  
\beq\label{uvcurve} 
U^2V+UV^2 +U+V-5UV+1=0, 
\eeq 
a plane cubic curve (total degree 3) which is also biquadratic (quadratic in both $U$ and $V$). 
The first orbit corresponds to the oval shown in red in Fig.\ref{s5orb}, whereas the second orbit lies outside the positive quadrant, moving around the three 
unbounded components of this curve, which can be seen in Fig.\ref{s5curvep}. 

\section{Heron triangles with two rational medians} 

\setcounter{equation}{0}

Buchholz found the first example of a Heron triangle with two rational medians: the $(73,51,26)$ triangle with area $420$ and  $k=\frac{35}{2}$, $\ell=\frac{97}{2}$, 
which had been overlooked by  Schubert in his work on parallelograms. 
After joining forces, Buchholz and  Rathbun conducted a systematic search for such triangles, using the following algorithm based on  
(\ref{abcparam}): fix the scale $\uptau=1$, enumerate pairs of rational numbers $(\theta,\phi)$, 
and for each pair use Heron's formula  (\ref{heron}) to check whether the area $\Delta$ is rational \cite{br1}. 
The first few triangles obtained from this search are shown in Table \ref{tab:table1}, where each triangle is represented by positive integers $(a,b,c)$ with $\gcd (a,b,c)=1$: 
their initial investigations suggested that there should be an infinite family 
of such  triangles (rows labelled with a positive integer $n$), 
together with an unknown number of sporadic triangles 
that do not fit into this family 
(rows labelled with asterisks).  

\begin{table}[h!]
  \begin{center}
    \caption{The smallest Heron triangles with two rational medians.}
    \label{tab:table1}
    \begin{tabular}{ | r|| r| r| r| r | r| r|} %
\hline
      $n$ & $a$  &  $b$ & $c$ & $k$ & $\ell$ & $\Delta$ \\
      \hline \hline
      1 & 73 & 51 & 26 & ${35}/{2}$ & ${97}/{2}$ & 420 \\
\hline
      2 & 626 &  875 & 291 & 572 &${433}/{2}$ & 55440 \\
\hline
      *    & 1241 & 4368 & 3673  & ${7975}/{2}$ & 1657 & 2042040 \\
\hline 
    **   &14384 & 14791 & 11257  & 11001 & 21177/2 & 75698280 \\
\hline
3 & 28779 & 13816 &  15155 & 3589/2 & 21937 & 23931600 \\ 
\hline
4 & 1823675& 185629& 1930456& 2048523/2& 3751059/2& 142334216640 \\
\hline
*** & 2288232 & 1976471 & 2025361 & 1641725 & 3843143/2 & 1877686881840 \\ 
\hline
**** & 22816608 & 20565641 & 19227017 & 
            16314487 & 36845705/2 & 185643608470320  \\ 
\hline 
5& 2442655864& 2396426547& 46263061& 1175099279& 2488886435/2& 2137147184560080 \\
\hline 
    \end{tabular}
  \end{center}
\end{table}

Heron triangles with two rational medians are associated with two different triples of Schubert parameters 
$(M_a,P_a,X_a)$, $(M_b,P_b,X_b)$, each corresponding to a particular set of angles $\al,\be,\gam$ adjacent to one of the medians $k$, $\ell$, respectively. 
These triples 
provide two different rational points on the Schubert surface (\ref{schub}), coupled by two constraints 
coming from the ratios of side lengths, as in (\ref{aboverc}). 
Remarkably, by considering the patterns of prime factors appearing in these rational numbers, Buchholz and Rathbun found conjectural formulae  
for a subset 
of these parameter triples 
in terms of the two Somos sequences (\ref{s5seq}), 
such as 
\beq\label{schubfac} 
M_a= -\frac{S_{n+1}S_{n+2}^2T_n}{S_nT_{n+1}T_{n+2}^2}, \qquad M_b=\frac{S_{n+1}S_{n+4}T_{n+1}T_{n+4}}{S_{n+2}S_{n+3}T_{n+2}T_{n+3}}, 
\eeq 
and analogous expressions for the other elements of each triple. When, for successive integers $n=1,2,3,\ldots$, they plotted the corresponding pairs $(\theta,\phi)$ 
found from  (\ref{thetaphi}) with a fixed choice of $\pm$ signs, they found them to lie on one of five algebraic curves ${\cal C}_j$, $1\leq j\leq5$, isomorphic to 
one another and repeating 
in a pattern with period 7, the simplest-looking curve being 
the biquadratic cubic 
\beq\label{curvec}
{\cal C}_4: \qquad 
\theta^2\phi-\theta\phi^2+\theta\phi+2\theta-2\phi-1=0. 
\eeq  

It was pointed out by Elkies 
that the sequences (\ref{s5seq}) can be written  using 
theta functions 
associated with the  
elliptic curve given by the 
equation 
\beq\label{ellipticc} 
E(\Q): \qquad y^2+xy=x^3+x^2-2x,
\eeq  
which has infinitely many rational points,\footnote{See  {\tt https://www.lmfdb.org/EllipticCurve/Q/102a1/}.} and is 
isomorphic (birationally equivalent) to ${\cal C}_4$. 
Indirectly, this led to a proof that 
every rational point $(\theta,\phi)$ on the genus one curve ${\cal C}_4$ given by (\ref{curvec}), 
subject to the constraints (\ref{thphcon}),  
corresponds 
to a Heron triangle with two rational medians \cite{br2}.

\begin{table}[h!]
  \begin{center}
    \caption{Prime factors of the semiperimeter, reduced side lengths  and area in the infinite family.} 
    \label{tab:table4}
\scalebox{0.85}{
    \begin{tabular}{ | r|| r| r| r| r | r|} %
\hline
      $n$ & $s$  &  $s-a$ & $s-b$ & $s-c$  & $\Delta$ \\
\hline 
    1 &  $3\cdot 5^2 $& 2 & $2^3\cdot 3$ & $7^2$  & $2^2\cdot3\cdot 5\cdot 7$ \\ 
      2 & $5\cdot 11^2$ & $3\cdot 7$ & $2\cdot 3^3\cdot 5$ &  $2^7\cdot 7$
& $2^4\cdot3^2\cdot 5\cdot 7\cdot 11$
 \\
	3 & $11\cdot 37^2$ & $2^3\cdot 5\cdot 7^3$& $3\cdot 5^3\cdot 7\cdot  11$ & $ 2^5\cdot 3$
& $2^4\cdot3\cdot 5^2\cdot 7^2\cdot 11\cdot 37$
\\
4 & $ 7\cdot 37 \cdot 83^2$ & $ 2^9\cdot 7 \cdot 11$ & $2^3\cdot 5 \cdot 11^3\cdot 37$ &  $3^4\cdot 5 \cdot 19^2$
& $2^6\cdot3^2\cdot 5\cdot 7\cdot 11^2\cdot19\cdot  37\cdot 83$
 \\  
5& $ 2^5\cdot 7^2\cdot 83\cdot 137^2$ & $2^3\cdot 3 \cdot 19 \cdot 37$ &  $11\cdot 37^3\cdot 83$ & $3\cdot 5^2 \cdot 11 \cdot 17^2 \cdot 19\cdot 23^2$
& $2^4\cdot3\cdot 5\cdot 7\cdot 11\cdot 17\cdot19\cdot 23\cdot  37^2\cdot 83\cdot 137$
  \\
\hline 
    \end{tabular}
}
  \end{center}
\end{table}

However, until very recently,  (\ref{schubfac}) and the explicit expressions for the  other  Schubert parameters remained conjectural.  
The key to progress  in \cite{honeheron} was to observe the elegant factorization pattern in the 
quantities appearing under the square root in Heron's formula, namely the semiperimeter $s$ and the reduced side lengths $s-a$, $s-b$, $s-c$  (see Table \ref{tab:table4}).  
It turns out that (up to an overall sign) each of these four quantities is given by a specific product of 
terms from the two Somos-5 sequences, leading to the following result.

\begin{thm}\label{main} 
For each integer $n\geq 1$, the terms in the pair of Somos-5 sequences $(S_n)$ and $(T_n)$ in (\ref{s5seq}) 
provide a Heron triangle with two rational medians, having integer side lengths given 
by  
$$ 
\bear{rcl} 
a & = & |S_{n+1}S_{n+2}^3S_{n+3}T_{n+2}+S_n^2S_{n+1}T_{n+3}T_{n+4}^2|, \\
b & = &  
| S_n^2S_{n+1}T_{n+3}T_{n+4}^2-T_{n+1}T_{n+2}^3T_{n+3}S_{n+2}|, \\ 
c & = &  
|T_{n+1}T_{n+2}^3T_{n+3}S_{n+2}-S_{n+1}S_{n+2}^3S_{n+3}T_{n+2}|, 
\eear 
$$  
with $\gcd(a,b,c)=1$, 
rational 
median lengths 
$$ 
\bear{rcl} 
k & = & {\scriptstyle \frac{1}{2} }
|S_{n+4}T_{n+4}(T_{n}T_{n+1}^2T_{n+2}-S_{n}S_{n+1}^2S_{n+2})|, \\
\ell & = &   {\scriptstyle \frac{1}{2} }
|S_nT_n(T_{n+2}T_{n+3}^2T_{n+4}-S_{n+2}S_{n+3}^2S_{n+4})|, 
\eear 
$$  
and integer area
$$\Delta = 
|S_nS_{n+1}S_{n+2}^2S_{n+3}S_{n+4}T_nT_{n+1}T_{n+2}^2T_{n+3}T_{n+4}|
.$$
\end{thm} 

\begin{table}[h!]
  \begin{center}
    \caption{Prime factors of the semiperimeter, reduced side lengths and area in the sporadic cases.}
    \label{tab:sporadicsabcbar}
\scalebox{0.8}{
    \begin{tabular}{ | r|| r| r| r| r | r |} %
\hline
      $n$ & $s$  &  $s-a$ & $s-b$ & $s-c$ & $\Delta$ \\
\hline 
    * &  $3\cdot 7\cdot 13\cdot 17 $& $2^3\cdot 5^2\cdot 17$ & $3\cdot 7\cdot 13$ &  $2^3\cdot 11^2$ & $2^3\cdot3\cdot 5\cdot 7\cdot 11\cdot 13\cdot 17$  \\ 
      ** & $2^3\cdot 7\cdot 19^2$ &  $2^3\cdot 3^6$  & $5^2\cdot 7\cdot 31$ 
&  $17^2\cdot 31$
& $2^3\cdot 3^3\cdot 5\cdot 7\cdot 17\cdot 19\cdot 31$ 
 \\
	*** & $2^3\cdot 3^2\cdot 11^2\cdot 19^2$ & $2^5\cdot 3^2\cdot 5^2\cdot 7\cdot 17$& $23^2\cdot 47^2$ & $ 7\cdot 17\cdot 97^2$ & $2^4\cdot3^2\cdot 5\cdot7\cdot 11\cdot 17 \cdot 19\cdot 23\cdot 47\cdot 97$
  \\
**** &   $17\cdot 23^2\cdot59^2$    
& $5^2\cdot 7^2\cdot 13^2\cdot 41$ 
& $2^4\cdot 3\cdot11^2\cdot 43^2$ 
& $2^4\cdot 3\cdot 17\cdot 19^2\cdot 41$
& $2^4\cdot 3\cdot 5\cdot 7\cdot 11\cdot 13\cdot 17\cdot 19 \cdot 23 \cdot 41\cdot 43 \cdot 59$
\\
\hline 
    \end{tabular}
}
  \end{center}
\end{table}

The curve (\ref{uvcurve}) is birationally equivalent to the curve ${\cal C}_4$ in (\ref{curvec}), hence also to the curve (\ref{ellipticc}). 
Its set of rational points is the union of two orbits of the map (\ref{qrts5}): the orbit associated with the sequence $(S_n)$, lying on the oval in the positive quadrant in 
Fig.\ref{s5curvep}; and the orbit associated with  $(T_n)$, which jumps around the other three quadrants in a pattern that repeats with period 7. 
Thus these two Somos-5 sequences completely encode the structure of this infinite family of Heron triangles with two rational medians. 
It is natural to wonder if any of the triangles in this family can have a third rational median $m$, but it has been proven that this is not the case \cite{ismail2}. 

Still, this leaves some big challenges to the reader: so far only four sporadic triangles have been found, which do not belong to the infinite family!
The prime factorization of each semiperimeter and the reduced lengths in Table \ref{tab:sporadicsabcbar} give tantalizing hints 
of further structure. Can you extend the search to find more sporadic examples, and fit them into one or more new infinite families, encoded by 
Somos (or other) sequences? Or can you show that these four are the only sporadic triangles, thereby proving that unicorns do not exist?

\end{document}